\documentclass[10pt]{article}
\usepackage{amsmath}
\usepackage{amssymb}
\usepackage{amsthm}
\usepackage{mathrsfs}
\usepackage[compress,noadjust]{cite}
\numberwithin{equation}{section}

\begin{document}
\title{Some estimates of intrinsic square functions on weighted Herz-type Hardy spaces}
\author{Hua Wang \footnote{E-mail address: wanghua@pku.edu.cn.}\\
\footnotesize{Department of Mathematics, Zhejiang University, Hangzhou 310027, China}}
\date{}
\maketitle

\begin{abstract}
In this paper, by using the atomic decomposition theory of weighted Herz-type Hardy spaces, we will obtain some strong type and weak type estimates for intrinsic square functions including the Lusin area function, Littlewood-Paley $\mathcal G$-function and $\mathcal G^*_\lambda$-function on these spaces.
\end{abstract}

\section{Introduction and main results}

Let ${\mathbb R}^{n+1}_+=\mathbb R^n\times(0,\infty)$ and $\varphi_t(x)=t^{-n}\varphi(x/t)$. The classical square function (Lusin area integral) is a familiar object. If $u(x,t)=P_t*f(x)$ is the Poisson integral of $f$, where $P_t(x)=\frac{c_nt}{(t^2+|x|^2)^{{(n+1)}/2}}$ denotes the Poisson kernel in ${\mathbb R}^{n+1}_+$. Then we define the classical square function (Lusin area integral) $S(f)$ by (see \cite{gundy} and \cite{stein})
\begin{equation*}
S(f)(x)=\bigg(\iint_{\Gamma(x)}\big|\nabla u(y,t)\big|^2t^{1-n}\,dydt\bigg)^{1/2},
\end{equation*}
where $\Gamma(x)$ denotes the usual cone of aperture one:
\begin{equation*}
\Gamma(x)=\big\{(y,t)\in{\mathbb R}^{n+1}_+:|x-y|<t\big\}
\end{equation*}
and
\begin{equation*}
\big|\nabla u(y,t)\big|^2=\left|\frac{\partial u}{\partial t}\right|^2+\sum_{j=1}^n\left|\frac{\partial u}{\partial y_j}\right|^2.
\end{equation*}
We can similarly define a cone of aperture $\gamma$ for any $\gamma>0$:
\begin{equation*}
\Gamma_\gamma(x)=\big\{(y,t)\in{\mathbb R}^{n+1}_+:|x-y|<\gamma t\big\},
\end{equation*}
and corresponding square function
\begin{equation*}
S_\gamma(f)(x)=\bigg(\iint_{\Gamma_\gamma(x)}\big|\nabla u(y,t)\big|^2t^{1-n}\,dydt\bigg)^{1/2}.
\end{equation*}
The Littlewood-Paley $g$-function (could be viewed as a ``zero-aperture" version of $S(f)$) and the $g^*_\lambda$-function (could be viewed as an ``infinite aperture" version of $S(f)$) are defined respectively by
(see, for example, \cite{muckenhoupt2} and \cite{stein2})
\begin{equation*}
g(f)(x)=\bigg(\int_0^\infty\big|\nabla u(x,t)\big|^2 t\,dt\bigg)^{1/2}
\end{equation*}
and
\begin{equation*}
g^*_\lambda(f)(x)=\left(\iint_{{\mathbb R}^{n+1}_+}\bigg(\frac t{t+|x-y|}\bigg)^{\lambda n}\big|\nabla u(y,t)\big|^2 t^{1-n}\,dydt\right)^{1/2}, \quad \lambda>1.
\end{equation*}

The modern (real-variable) variant of $S_\gamma(f)$ can be defined in the following way (here we drop the subscript $\gamma$ if $\gamma=1$). Let $\psi\in C^\infty(\mathbb R^n)$ be real, radial, have support contained in $\{x:|x|\le1\}$, and $\int_{\mathbb R^n}\psi(x)\,dx=0$. The continuous square function $S_{\psi,\gamma}(f)$ is defined by (see, for instance, \cite{chang} and \cite{chanillo})
\begin{equation*}
S_{\psi,\gamma}(f)(x)=\bigg(\iint_{\Gamma_\gamma(x)}\big|f*\psi_t(y)\big|^2\frac{dydt}{t^{n+1}}\bigg)^{1/2}.
\end{equation*}

In 2007, Wilson \cite{wilson1} introduced a new square function called intrinsic square function which is universal in a sense (see also \cite{wilson2}). This function is independent of any particular kernel $\psi$, and it dominates pointwise all the above-defined square functions. On the other hand, it is not essentially larger than any particular $S_{\psi,\gamma}(f)$. For $0<\beta\le1$, let ${\mathcal C}_\beta$ be the family of functions $\varphi$ defined on $\mathbb R^n$ such that $\varphi$ has support containing in $\{x\in\mathbb R^n: |x|\le1\}$, $\int_{\mathbb R^n}\varphi(x)\,dx=0$, and for all $x, x'\in \mathbb R^n$,
\begin{equation*}
|\varphi(x)-\varphi(x')|\le|x-x'|^\beta.
\end{equation*}
For $(y,t)\in {\mathbb R}^{n+1}_+$ and $f\in L^1_{{loc}}(\mathbb R^n)$, we set
\begin{equation}
A_\beta(f)(y,t)=\sup_{\varphi\in{\mathcal C}_\beta}\big|f*\varphi_t(y)\big|=\sup_{\varphi\in{\mathcal C}_\beta}\bigg|\int_{\mathbb R^n}\varphi_t(y-z)f(z)\,dz\bigg|.
\end{equation}
Then we define the intrinsic square function of $f$ (of order $\beta$) by the formula
\begin{equation}
\mathcal S_\beta(f)(x)=\left(\iint_{\Gamma(x)}\Big(A_\beta(f)(y,t)\Big)^2\frac{dydt}{t^{n+1}}\right)^{1/2}.
\end{equation}
We can also define varying-aperture versions of $\mathcal S_\beta(f)$ by the formula
\begin{equation}
\mathcal S_{\beta,\gamma}(f)(x)=\left(\iint_{\Gamma_\gamma(x)}\Big(A_\beta(f)(y,t)\Big)^2\frac{dydt}{t^{n+1}}\right)^{1/2}.
\end{equation}
The intrinsic Littlewood-Paley $\mathcal G$-function and the intrinsic $\mathcal G^*_\lambda$-function will be given respectively by
\begin{equation}
\mathcal G_\beta(f)(x)=\left(\int_0^\infty\Big(A_\beta(f)(x,t)\Big)^2\frac{dt}{t}\right)^{1/2}
\end{equation}
and
\begin{equation}
\mathcal G^*_{\lambda,\beta}(f)(x)=\left(\iint_{{\mathbb R}^{n+1}_+}\left(\frac t{t+|x-y|}\right)^{\lambda n}\Big(A_\beta(f)(y,t)\Big)^2\frac{dydt}{t^{n+1}}\right)^{1/2}, \lambda>1.
\end{equation}

In \cite{wilson2}, Wilson showed the following weighted $L^p$ boundedness of the intrinsic square functions.

\newtheorem*{thma}{Theorem A}

\begin{thma}
Let $0<\beta\le1$, $1<p<\infty$ and $w\in A_p (\mbox{Muckenhoupt weight class})$. Then there exists a constant $C>0$ independent of $f$ such that
\begin{equation*}
\big\|\mathcal S_\beta(f)\big\|_{L^p_w}\le C \|f\|_{L^p_w}.
\end{equation*}
\end{thma}
Moreover, in \cite{lerner}, Lerner obtained sharp $L^p_w$ norm inequalities for the intrinsic square functions in terms of the $A_p$ characteristic constant of $w$ for all $1<p<\infty$. For further discussions about the boundedness of intrinsic square functions on various function spaces, we refer the readers to \cite{huang,wang2,wang3,wang4,wang5,wang6}.

The aim of this paper is to discuss the boundedness properties of intrinsic square functions on the homogeneous (non-homogeneous) weighted Herz-type Hardy spaces (see Section 2 below for the definitions). Moreover, at the endpoint case, we will obtain their weak type estimates. Our main results are stated as follows.

\newtheorem{theorem}{Theorem}[section]

\begin{theorem}
Let $w_1$, $w_2\in A_1$, $0<p<\infty$, $1<q<\infty$, $0<\beta\le1$ and $n(1-1/q)\le\alpha<n(1-1/q)+\beta$. Then $\mathcal S_\beta$ is bounded from $H\dot K^{\alpha,p}_q(w_1,w_2)$ $(HK^{\alpha,p}_q(w_1,w_2))$ into $\dot K^{\alpha,p}_q(w_1,w_2)$ $(K^{\alpha,p}_q(w_1,w_2))$.
\end{theorem}

\begin{theorem}
Let $w_1$, $w_2\in A_1$, $0<p\le1$, $1<q<\infty$, $0<\beta<1$ and $\alpha=n(1-1/q)+\beta$. Then $\mathcal S_\beta$ is bounded from $H\dot K^{\alpha,p}_q(w_1,w_2)$ $(HK^{\alpha,p}_q(w_1,w_2))$ into $W\dot K^{\alpha,p}_q(w_1,w_2)$ $(WK^{\alpha,p}_q(w_1,w_2))$.
\end{theorem}

\begin{theorem}
Let $w_1$, $w_2\in A_1$, $0<p<\infty$, $1<q<\infty$, $0<\beta\le1$ and $n(1-1/q)\le\alpha<n(1-1/q)+\beta$. Suppose that $\lambda>3+(2\beta)/n$, then $\mathcal G^*_{\lambda,\beta}$ is bounded from $H\dot K^{\alpha,p}_q(w_1,w_2)$ $(HK^{\alpha,p}_q(w_1,w_2))$ into $\dot K^{\alpha,p}_q(w_1,w_2)$ $(K^{\alpha,p}_q(w_1,w_2))$.
\end{theorem}

\begin{theorem}
Let $w_1$, $w_2\in A_1$, $0<p\le1$, $1<q<\infty$, $0<\beta<1$ and $\alpha=n(1-1/q)+\beta$. Suppose that $\lambda>3+(2\beta)/n$, then $\mathcal G^*_{\lambda,\beta}$ is bounded from $H\dot K^{\alpha,p}_q(w_1,w_2)$ $(HK^{\alpha,p}_q(w_1,w_2))$ into $W\dot K^{\alpha,p}_q(w_1,w_2)$ $(WK^{\alpha,p}_q(w_1,w_2))$.
\end{theorem}

In \cite{wilson1}, Wilson also showed that for any $0<\beta\le1$, the functions $\mathcal S_\beta(f)(x)$ and $\mathcal G_\beta(f)(x)$ are pointwise comparable, with comparability constants depending only on $\beta$ and $n$. Thus, as a direct consequence of Theorems 1.1 and 1.2, we obtain the following:

\newtheorem{corollary}[theorem]{Corollary}

\begin{corollary}
Let $w_1$, $w_2\in A_1$, $0<p<\infty$, $1<q<\infty$, $0<\beta\le1$ and $n(1-1/q)\le\alpha<n(1-1/q)+\beta$. Then $\mathcal G_\beta$ is bounded from $H\dot K^{\alpha,p}_q(w_1,w_2)$ $(HK^{\alpha,p}_q(w_1,w_2))$ into $\dot K^{\alpha,p}_q(w_1,w_2)$ $(K^{\alpha,p}_q(w_1,w_2))$.
\end{corollary}

\begin{corollary}
Let $w_1$, $w_2\in A_1$, $0<p\le1$, $1<q<\infty$, $0<\beta<1$ and $\alpha=n(1-1/q)+\beta$. Then $\mathcal G_\beta$ is bounded from $H\dot K^{\alpha,p}_q(w_1,w_2)$ $(HK^{\alpha,p}_q(w_1,w_2))$ into $W\dot K^{\alpha,p}_q(w_1,w_2)$ $(WK^{\alpha,p}_q(w_1,w_2))$.
\end{corollary}

\section{Notations and preliminaries}

\subsection{$A_p$ weights}

Let us first recall some standard definitions and notations. The classical $A_p$ weight
theory was first introduced by Muckenhoupt in the study of weighted
$L^p$ boundedness of Hardy-Littlewood maximal functions in \cite{muchenhoupt1}. A weight $w$ is a locally integrable function on $\mathbb R^n$ which takes values in $(0,\infty)$ almost everywhere. $B=B(x_0,R)$ denotes the ball with the center $x_0$ and radius $R$. Given a ball $B$ and $\lambda>0$, $\lambda B$ stands for the ball concentric with $B$ whose radius is $\lambda$ times as long. For a given weight function $w$ and a measurable set $E$, we also denote the Lebesgue measure of $E$ by $|E|$ and the weighted measure of $E$ by $w(E)$, where $w(E)=\int_E w(x)\,dx$. We say that $w$ is in the Muckenhoupt class $A_p$ with $1<p<\infty$, if
\begin{equation*}
\left(\frac1{|B|}\int_B w(x)\,dx\right)\left(\frac1{|B|}\int_B w(x)^{-1/{(p-1)}}\,dx\right)^{p-1}\le C, \quad\mbox{for every ball}\; B\subseteq \mathbb R^n,
\end{equation*}
where $C$ is a positive constant which is independent of the choice of $B$. For the case $p=1$, $w\in A_1$, if
\begin{equation*}
\frac1{|B|}\int_B w(x)\,dx\le C\cdot\underset{x\in B}{\mbox{ess\,inf}}\,w(x),
\quad\mbox{for every ball}\;B\subseteq\mathbb R^n.
\end{equation*}
The smallest value of $C$ such that the above inequality holds is called the $A_1$ characteristic constant of $w$ and denoted by $[w]_{A_1}$. A weight function $w$ is said to belong to the reverse H\"{o}lder class $RH_r$, if there exist two constants $r>1$ and $C>0$ such that the following reverse H\"{o}lder inequality holds
\begin{equation*}
\left(\frac{1}{|B|}\int_B w(x)^r\,dx\right)^{1/r}\le C\left(\frac{1}{|B|}\int_B w(x)\,dx\right),
\quad\mbox{for every ball}\; B\subseteq \mathbb R^n.
\end{equation*}

It is well known that if $w\in A_p$ with $1<p<\infty$, then $w\in A_r$ for all $r>p$, and $w\in A_q$ for some $1<q<p$. Moreover, if $w\in A_p$ with $1\le p<\infty$, then there exists $r>1$ such that $w\in RH_r$.

We state the following results that we will use frequently in the sequel.

\newtheorem{lemma}[theorem]{Lemma}

\begin{lemma}[\cite{garcia2}]
Let $w\in A_1$. Then, for any ball $B$, there exists an absolute constant $C>0$ such that
\begin{equation*}
w(2B)\le C\,w(B).
\end{equation*}
More precisely, for any $\lambda>1$, we have
\begin{equation*}
w(\lambda B)\le [w]_{A_1}\cdot\lambda^{n}w(B).
\end{equation*}
\end{lemma}

\begin{lemma}[\cite{garcia2,gundy}]
Let $w\in A_1\cap RH_r$, $r>1$. Then there exist two constants $C_1$, $C_2>0$ such that
\begin{equation*}
C_1\left(\frac{|E|}{|B|}\right)\le\frac{w(E)}{w(B)}\le C_2\left(\frac{|E|}{|B|}\right)^{(r-1)/r}
\end{equation*}
for any measurable subset $E$ of a ball $B$.
\end{lemma}

\subsection{Weighted Herz-type Hardy spaces}

Next we shall give the definitions of the weighted Herz space, weighted weak Herz space and weighted Herz-type Hardy space. In 1964, Beurling \cite{beurling} first introduced some fundamental form of Herz spaces to study convolution algebras. Later Herz \cite{herz} gave versions of the spaces defined below in a slightly different setting. Since then, the theory of Herz spaces has been significantly developed, and these spaces have turned out to be quite useful in harmonic analysis. For instance, they were used by Baernstein and Sawyer \cite{baernstein} to characterize the multipliers on the classical Hardy spaces, and used by Lu and Yang \cite{lu6,lu7} in the study of partial differential equations. The weighted version of Herz spaces was also introduced and investigated in \cite{komori,lu1,lu2,lu6,tang}.

On the other hand, a theory of Hardy spaces associated with Herz spaces has been developed in \cite{garcia1,lu4}. These new Herz-type Hardy spaces may be regarded as the local version at the origin of the classical Hardy spaces $H^p(\mathbb R^n)$ and are good substitutes for $H^p(\mathbb R^n)$ when we study the boundedness of non-translation invariant operators (see \cite{grafakos,lu5,lu9}). For the weighted case, in 1995, Lu and Yang \cite{lu3,lu8} introduced the following weighted Herz-type Hardy spaces $H\dot K^{\alpha,p}_q(w_1,w_2)$
($H K^{\alpha,p}_q(w_1,w_2)$) and established their central atomic decompositions. For further details about the properties and boundedness of some operators on weighted Herz-type Hardy spaces, we refer the readers to \cite{lee1,lee2,li,lu10,wang1} and the references therein.

Let $B_k=\{x\in\mathbb R^n:|x|\le 2^k\}$ and $C_k=B_k\backslash B_{k-1}$ for $k\in\mathbb Z$. Denote $\chi_k=\chi_{_{C_k}}$ for $k\in\mathbb Z$, $\widetilde\chi_k=\chi_k$ if $k\in\mathbb N$ and $\widetilde\chi_0=\chi_{_{B_0}}$, where $\chi_{E}$ is the characteristic function of a set $E$. For any given weight function $w$ on $\mathbb R^n$ and $0<q<\infty$, we denote by $L^q_w(\mathbb R^n)$ the space of all functions $f$ satisfying
\begin{equation}
\|f\|_{L^q_w}=\left(\int_{\mathbb R^n}|f(x)|^qw(x)\,dx\right)^{1/q}<\infty.
\end{equation}

\newtheorem{def1}[theorem]{Definition}

\begin{def1}[\cite{lu2}]
Let $\alpha\in\mathbb R$, $0<p, q<\infty$ and $w_1$, $w_2$ be two weight functions on $\mathbb R^n$.

$(a)$ The homogeneous weighted Herz space $\dot K^{\alpha,p}_q(w_1,w_2)$ is defined by
\begin{equation*}
\dot K^{\alpha,p}_q(w_1,w_2)=\Big\{f\in L^q_{loc}(\mathbb R^n\backslash\{0\},w_2):\big\|f\big\|_{\dot K^{\alpha,p}_q(w_1,w_2)}<\infty\Big\},
\end{equation*}
where
\begin{equation}
\big\|f\big\|_{\dot K^{\alpha,p}_q(w_1,w_2)}=\left(\sum_{k\in\mathbb Z}\big[w_1(B_k)\big]^{{\alpha p}/n}\big\|f\chi_k\big\|_{L^q_{w_2}}^p\right)^{1/p}.
\end{equation}

$(b)$ The non-homogeneous weighted Herz space $K^{\alpha,p}_q(w_1,w_2)$ is defined by
\begin{equation*}
K^{\alpha,p}_q(w_1,w_2)=\Big\{f\in L^q_{loc}(\mathbb R^n,w_2):\big\|f\big\|_{K^{\alpha,p}_q(w_1,w_2)}<\infty\Big\},
\end{equation*}
where
\begin{equation}
\big\|f\big\|_{K^{\alpha,p}_q(w_1,w_2)}=\left(\sum_{k=0}^\infty\big[w_1(B_k)\big]^{{\alpha p}/n}\big\|f\widetilde{\chi}_k\big\|_{L^q_{w_2}}^p\right)^{1/p}.
\end{equation}
\end{def1}

For any $k\in\mathbb Z$, $\lambda>0$ and any measurable function $f$ on $\mathbb R^n$, we set $E_k(\lambda,f)=\{x\in C_k:|f(x)|>\lambda\}$. Let $\widetilde E_k(\lambda,f)=E_k(\lambda,f)$ for $k\in\mathbb N$ and $\widetilde E_0(\lambda,f)=\{x\in B(0,1):|f(x)|>\lambda\}$.

\begin{def1}[\cite{lu1}]
Let $\alpha\in\mathbb R$, $0<p, q<\infty$ and $w_1$, $w_2$ be two weight functions on $\mathbb R^n$.

$(c)$ A measurable function $f(x)$ on $\mathbb R^n$ is said to belong to the homogeneous weighted weak Herz space $W\dot K^{\alpha,p}_q(w_1,w_2)$ if
\begin{equation}
\big\|f\big\|_{W\dot K^{\alpha,p}_q(w_1,w_2)}=\sup_{\lambda>0}\lambda\left(\sum_{k\in\mathbb Z}\big[w_1(B_k)\big]^{{\alpha p}/n}\big[w_2(E_k(\lambda,f))\big]^{p/q}\right)^{1/p}<\infty.
\end{equation}

$(d)$ A measurable function $f(x)$ on $\mathbb R^n$ is said to belong to the non-homogeneous weighted weak Herz space $W K^{\alpha,p}_q(w_1,w_2)$ if
\begin{equation}
\big\|f\big\|_{W K^{\alpha,p}_q(w_1,w_2)}=\sup_{\lambda>0}\lambda\left(\sum_{k=0}^\infty\big[w_1(B_k)\big]^{{\alpha p}/n}\big[w_2(\widetilde E_k(\lambda,f))\big]^{p/q}\right)^{1/p}<\infty.
\end{equation}
\end{def1}

Let $\mathscr S(\mathbb R^n)$ be the class of Schwartz functions and let $\mathscr S'(\mathbb R^n)$ be its dual space. For any given $f\in\mathscr S'(\mathbb R^n)$, then the grand maximal function of $f$ is defined by
\begin{equation*}
G(f)(x)=\sup_{\varphi\in{\mathscr A_N}}\sup_{|y-x|<t}\big|\varphi_t*f(y)\big|,
\end{equation*}
where $\mathscr A_N=\Big\{\varphi\in\mathscr S(\mathbb R^n):\sup_{|\alpha|,|\beta|\le N}\big|x^\alpha D^\beta\varphi(x)\big|\le1\Big\}$ and $N\in\mathbb N$ is sufficiently large.

\begin{def1}[\cite{lu3}]
Let $0<\alpha<\infty$, $0<p<\infty$, $1<q<\infty$ and $w_1$, $w_2$ be two weight functions on $\mathbb R^n$.

$(e)$ The homogeneous weighted Herz-type Hardy space $H\dot K^{\alpha,p}_q(w_1,w_2)$ associated with the space $\dot K^{\alpha,p}_q(w_1,w_2)$ is defined by
\begin{equation*}
H\dot K^{\alpha,p}_q(w_1,w_2)=\Big\{f\in\mathscr S'(\mathbb R^n):G(f)\in\dot K^{\alpha,p}_q(w_1,w_2)\Big\}
\end{equation*}
and we define $\big\|f\big\|_{H\dot K^{\alpha,p}_q(w_1,w_2)}=\big\|G(f)\big\|_{\dot K^{\alpha,p}_q(w_1,w_2)}$.

$(f)$ The non-homogeneous weighted Herz-type Hardy space $HK^{\alpha,p}_q(w_1,w_2)$ associated with the space $K^{\alpha,p}_q(w_1,w_2)$ is defined by
\begin{equation*}
HK^{\alpha,p}_q(w_1,w_2)=\Big\{f\in\mathscr S'(\mathbb R^n):G(f)\in K^{\alpha,p}_q(w_1,w_2)\Big\}
\end{equation*}
and we define $\big\|f\big\|_{HK^{\alpha,p}_q(w_1,w_2)}=\big\|G(f)\big\|_{K^{\alpha,p}_q(w_1,w_2)}$.
\end{def1}

In this article, we will use Lu and Yang's central atomic decomposition theory for weighted Herz-type Hardy spaces in \cite{lu3,lu8} (see also \cite{lu10}). We characterize weighted Herz-type Hardy spaces in terms of central atoms in the following way.

\begin{def1}[\cite{lu3}]
Let $1<q<\infty$, $n(1-1/q)\le\alpha<\infty$ and $s\ge[\alpha+n(1/q-1)]$.
\\
$(i)$ A function $a(x)$ on $\mathbb R^n$ is said to be a central $(\alpha,q,s)$-atom with respect to $(w_1,w_2)$ $($or a central $(\alpha,q,s;w_1,w_2)$-atom$)$, if it satisfies

$(a)$ supp\,$a\subseteq B(0,R)=\{x\in\mathbb R^n:|x|\le R\}, \;R>0;$

$(b)\;\big\|a\big\|_{L^q_{w_2}}\le \big[w_1(B(0,R))\big]^{-\alpha/n};$

$(c)\;\int_{\mathbb R^n}a(x)x^\gamma\,dx=0, \;\mbox{for every multi-index}\;\gamma\;\mbox{with}\;|\gamma|\le s$.
\\
$(ii)$ A function $a(x)$ on $\mathbb R^n$ is said to be a central $(\alpha,q,s)$-atom of restricted type with respect to $(w_1,w_2)$ $($or a central $(\alpha,q,s;w_1,w_2)$-atom of restricted type$)$, if it satisfies the conditions $(b)$, $(c)$ above and

$(a')$ supp\,$a\subseteq B(0,R),$ for some $R>1$.
\end{def1}

\begin{theorem}[\cite{lu3}]
Let $w_1$, $w_2\in A_1$, $0<p<\infty$, $1<q<\infty$, $n(1-1/q)\le\alpha<\infty$ and $s\ge[\alpha+n(1/q-1)]$. Then we have

$(i)$ $f\in H\dot K^{\alpha,p}_q(w_1,w_2)$ if and only if
\begin{equation*}
f(x)=\sum_{j\in\mathbb Z}\lambda_ja_j(x),\quad \mbox{in the sense of}\;\,\mathscr S'(\mathbb R^n),
\end{equation*}
where $\sum_{j\in\mathbb Z}|\lambda_j|^p<\infty$, each $a_j$ is a central $(\alpha,q,s;w_1,w_2)$-atom with supp\,$a_j\subseteq B_j=B(0,2^j)$. Moreover,
\begin{equation*}
\big\|f\big\|_{H\dot K^{\alpha,p}_q(w_1,w_2)}\approx\inf\bigg(\sum_{j\in\mathbb Z}|\lambda_j|^p\bigg)^{1/p},
\end{equation*}
where the infimum is taken over all the above decompositions of $f$.

$(ii)$ $f\in HK^{\alpha,p}_q(w_1,w_2)$ if and only if
\begin{equation*}
f(x)=\sum_{j=0}^\infty\lambda_ja_j(x),\quad \mbox{in the sense of}\;\,\mathscr S'(\mathbb R^n),
\end{equation*}
where $\sum_{j=0}^\infty|\lambda_j|^p<\infty$, each $a_j$ is a central $(\alpha,q,s;w_1,w_2)$-atom of restricted type with supp\,$a_j\subseteq B_j=B(0,2^j)$. Moreover,
\begin{equation*}
\big\|f\big\|_{HK^{\alpha,p}_q(w_1,w_2)}\approx\inf\bigg(\sum_{j=0}^\infty|\lambda_j|^p\bigg)^{1/p},
\end{equation*}
where the infimum is taken over all the above decompositions of $f$.
\end{theorem}

Throughout this article, we will use $C$ to denote a positive constant, which is independent of the main parameters and not necessarily the same at each occurrence.

\section{Proofs of Theorems 1.1 and 1.2}

\begin{proof}[Proof of Theorem 1.1]
First we note that the assumptions $n(1-1/q)\le\alpha<n(1-1/q)+\beta$ and $0<\beta\le1$ imply that $N=[\alpha+n(1/q-1)]=0$. For any central $(\alpha,q,0;w_1,w_2)$-atom $a$ with supp\,$a\subseteq B_{\ell}=B(0,2^\ell)$, $\ell\in\mathbb Z$, we are going to show that $\big\|\mathcal S_\beta(a)\big\|_{\dot K^{\alpha,p}_q(w_1,w_2)}\le C$ for the case of $0<p\le1$, where $C>0$ is a universal constant independent of the choice of $a$. Write
\begin{equation*}
\begin{split}
\big\|\mathcal S_\beta(a)\big\|^p_{\dot K^{\alpha,p}_q(w_1,w_2)}&=
\sum_{k\in\mathbb Z}\big[w_1(B_k)\big]^{{\alpha p}/n}\big\|\mathcal S_\beta(a)\chi_k\big\|^p_{L^q_{w_2}}\\
&=\sum_{k=-\infty}^{\ell+1}\big[w_1(B_k)\big]^{{\alpha p}/n}\big\|\mathcal S_\beta(a)\chi_k\big\|^p_{L^q_{w_2}}
+\sum_{k=\ell+2}^{\infty}\big[w_1(B_k)\big]^{{\alpha p}/n}\big\|\mathcal S_\beta(a)\chi_k\big\|^p_{L^q_{w_2}}\\
&=I_1+I_2.
\end{split}
\end{equation*}
Since $w_2\in A_1$, then $w_2\in A_q$ for any $1<q<\infty$. It follows from Theorem A that
\begin{equation*}
I_1\le\sum_{k=-\infty}^{\ell+1}\big[w_1(B_k)\big]^{{\alpha p}/n}\big\|\mathcal S_\beta(a)\big\|^p_{L^q_{w_2}}
\le C\sum_{k=-\infty}^{\ell+1}\big[w_1(B_k)\big]^{{\alpha p}/n}\big\|a\big\|^p_{L^q_{w_2}}.
\end{equation*}
Since $w_1\in A_1$, then we know that $w\in RH_r$ for some $r>1$. When $k\le \ell+1$, then $B_k\subseteq B_{\ell+1}$. Consequently, by Lemma 2.2, we have
\begin{equation}
\frac{w_1(B_k)}{w_1(B_{\ell+1})}\le C\cdot\left(\frac{|B_k|}{|B_{\ell+1}|}\right)^{\delta},
\end{equation}
where $\delta=(r-1)/r>0$. Thus, by using the size condition of central atom $a$ and (3.1), we obtain
\begin{equation*}
\begin{split}
I_1&\le C\sum_{k=-\infty}^{\ell+1} 2^{(k-\ell-1)\alpha\delta p}\\
&=C\sum_{k=-\infty}^{0}2^{k\alpha\delta p}\\
&\le C.
\end{split}
\end{equation*}
To estimate the other term $I_2$, we first claim that for any $(y,t)\in{\mathbb R}^{n+1}_+$, the following inequality holds:
\begin{equation}
A_\beta(a)(y,t)\le C\cdot\frac{2^{\ell(n+\beta)}}{t^{n+\beta}}\big[w_1(B_\ell)\big]^{-\alpha/n}\big[w_2(B_\ell)\big]^{-1/q}.
\end{equation}
In fact, for any $\varphi\in{\mathcal C}_\beta$ with $0<\beta\le1$, by the vanishing moment condition of central atom $a$, we have
\begin{align}
\big|a*\varphi_t(y)\big|&=\left|\int_{B_\ell}\Big[\varphi_t(y-z)-\varphi_t(y)\Big]a(z)\,dz\right|\notag\\
&\le\int_{B_\ell}\frac{|z|^\beta}{t^{n+\beta}}\big|a(z)\big|\,dz\notag\\
&\le\frac{2^{\beta\ell}}{t^{n+\beta}}\int_{B_\ell}|a(z)|\,dz.
\end{align}
Denote the conjugate exponent of $q>1$ by $q'=q/(q-1)$. Using H\"older's inequality, $A_q$ condition and the size condition of central atom $a$, we can get
\begin{align}
\int_{B_\ell}|a(z)|\,dz&\le\bigg(\int_{B_\ell}\big|a(z)\big|^qw_2(z)\,dz\bigg)^{1/q}
\bigg(\int_{B_\ell} w_2(z)^{-{q'}/q}\,dz\bigg)^{1/{q'}}\notag\\
&\le C\cdot\big\|a\big\|_{L^q_{w_2}}\big|B_\ell\big|\big[w_2(B_\ell)\big]^{-1/q}\notag\\
&\le C\cdot\big|B_\ell\big|\big[w_1(B_\ell)\big]^{-\alpha/n}\big[w_2(B_\ell)\big]^{-1/q}.
\end{align}
Substituting the above inequality (3.4) into (3.3) and then taking the supremum over all functions $\varphi\in\mathcal C_\beta$, we obtain the desired inequality (3.2). Observe that if $x\in C_k=B_k\backslash B_{k-1}$, $k\ge\ell+2$ and $z\in B_\ell$, then we have $|z|\le\frac12|x|$. We also note that supp\,$\varphi\subseteq\{x\in\mathbb R^n:|x|\le1\}$, then for any $z\in B_\ell$, $(y,t)\in\Gamma(x)$ and $x\in C_k$ with $k\ge\ell+2$, we can deduce that
\begin{equation*}
2t>|x-y|+|y-z|\ge|x-z|\ge|x|-|z|\ge\frac{|x|}{2}.
\end{equation*}
Hence, for any $x\in C_k=B_k\backslash B_{k-1}$ with $k\ge\ell+2$, by using the inequality (3.2), we obtain
\begin{align}
\big|\mathcal S_\beta(a)(x)\big|&\le C\bigg(2^{\ell(n+\beta)}\big[w_1(B_\ell)\big]^{-\alpha/n}\big[w_2(B_\ell)\big]^{-1/q}\bigg)
\left(\int_{\frac{|x|}{4}}^\infty\int_{|y-x|<t}\frac{dydt}{t^{2n+2\beta+n+1}}\right)^{1/2}\notag\\
&\le C\bigg(2^{\ell(n+\beta)}\big[w_1(B_\ell)\big]^{-\alpha/n}\big[w_2(B_\ell)\big]^{-1/q}\bigg)
\left(\int_{\frac{|x|}{4}}^\infty\frac{dt}{t^{2n+2\beta+1}}\right)^{1/2}\notag\\
&\le C\cdot2^{\ell(n+\beta)}\big[w_1(B_\ell)\big]^{-\alpha/n}\big[w_2(B_\ell)\big]^{-1/q}
\cdot\frac{1}{|x|^{n+\beta}}.
\end{align}
Substituting the above inequality (3.5) into the term $I_2$, we can see that
\begin{equation*}
\begin{split}
I_2=&\sum_{k=\ell+2}^{\infty}\big[w_1(B_k)\big]^{{\alpha p}/n}\bigg(\int_{2^{k-1}<|x|\le 2^k}\big|\mathcal S_\beta(a)(x)\big|^qw_2(x)\,dx\bigg)^{p/q}\\
\le& C\sum_{k=\ell+2}^{\infty}\big[w_1(B_k)\big]^{{\alpha p}/n}\bigg(2^{\ell(n+\beta)}\big[w_1(B_\ell)\big]^{-\alpha/n}\big[w_2(B_\ell)\big]^{-1/q}\bigg)^p\\
\end{split}
\end{equation*}
\begin{equation*}
\begin{split}
&\times\bigg(\int_{2^{k-1}<|x|\le 2^k}\frac{w_2(x)}{|x|^{q(n+\beta)}}\,dx\bigg)^{p/q}\\
\le& C\sum_{k=\ell+2}^{\infty}\bigg(\frac{2^{\ell p(n+\beta)}}{2^{kp(n+\beta)}}\bigg
)\left(\frac{w_1(B_k)}{w_1(B_\ell)}\right)^{{\alpha p}/n}\left(\frac{w_2(B_k)}{w_2(B_\ell)}\right)^{p/q}.
\end{split}
\end{equation*}
In this case, when $k\ge\ell+2$, then we have $B_k\supseteq B_{\ell+2}\supseteq B_{\ell}$. Since $w_1,w_2\in A_1$, then by using Lemma 2.2 again, we can get
\begin{equation}
\frac{w_i(B_k)}{w_i(B_{\ell})}\le C\cdot \frac{|B_k|}{|B_{\ell}|}, \quad \mbox{for}\; i=1 \;\mbox{and}\;2.
\end{equation}
Hence, from the above inequality (3.6), it follows that
\begin{equation*}
\begin{split}
I_2&\le C\sum_{k=\ell+2}^{\infty}\left(\frac{2^{\ell p(n+\beta)}}{2^{kp(n+\beta)}}\right)
\left(\frac{2^{kn}}{2^{\ell n}}\right)^{{\alpha p}/n}\left(\frac{2^{kn}}{2^{\ell n}}\right)^{p/q}\notag\\
&=C\sum_{k=2}^{\infty}\left(\frac{1}{2^k}\right)^{p(n+\beta)-\alpha p-{np}/q}\notag\\
&\le C,
\end{split}
\end{equation*}
where the last series is convergent since $\alpha<n(1-1/q)+\beta$. Combining the above estimates for $I_1$ and $I_2$, we get the desired result.

We are now in a position to give the proof of Theorem 1.1 for the case $0<p\le1$. For every $f\in H\dot K^{\alpha,p}_q(w_1,w_2)$, then by Theorem 2.7, we have the decomposition $f=\sum_{\ell\in\mathbb Z}\lambda_\ell a_\ell$, where $\sum_{\ell\in\mathbb Z}|\lambda_\ell|^p<\infty$ and each $a_\ell$ is a central $(\alpha,q,0;w_1,w_2)$-atom with supp\,$a_\ell\subseteq B_\ell=B(0,2^{\ell})$. Therefore
\begin{equation*}
\begin{split}
\big\|\mathcal S_\beta(f)\big\|^p_{\dot K^{\alpha,p}_q(w_1,w_2)}
&\le C\sum_{k\in\mathbb Z}\big[w_1(B_k)\big]^{{\alpha p}/n}\Bigg(\sum_{\ell\in\mathbb Z}\big|\lambda_\ell\big|\big\|\mathcal S_\beta(a_\ell)\chi_k\big\|_{L^q_{w_2}}\Bigg)^p\\
&\le C\sum_{k\in\mathbb Z}\big[w_1(B_k)\big]^{{\alpha p}/n}\Bigg(\sum_{\ell\in\mathbb Z}\big|\lambda_\ell\big|^p\big\|\mathcal S_\beta(a_\ell)\chi_k\big\|^p_{L^q_{w_2}}\Bigg)\\
&\le C\sum_{\ell\in\mathbb Z}\big|\lambda_\ell\big|^p\\
&\le C\big\|f\big\|^p_{H\dot K^{\alpha,p}_q(w_1,w_2)}.
\end{split}
\end{equation*}

We now consider the case $1<p<\infty$. As above, we write
\begin{equation*}
\begin{split}
\big\|\mathcal S_\beta(f)\big\|^p_{\dot K^{\alpha,p}_q(w_1,w_2)}
\le&\,\sum_{k\in\mathbb Z}\big[w_1(B_k)\big]^{{\alpha p}/n}
\Bigg(\sum_{\ell=k-1}^{\infty}\big|\lambda_\ell\big|
\big\|\mathcal S_\beta(a_\ell)\chi_k\big\|_{L^q_{w_2}}\Bigg)^p\\
&+\sum_{k\in\mathbb Z}\big[w_1(B_k)\big]^{{\alpha p}/n}
\Bigg(\sum_{\ell=-\infty}^{k-2}\big|\lambda_\ell\big|
\big\|\mathcal S_\beta(a_\ell)\chi_k\big\|_{L^q_{w_2}}\Bigg)^p\\
=&\,I_1'+I_2'.
\end{split}
\end{equation*}
Let us first deal with $I_1'$. Applying H\"older's inequality, Theorem A and the size condition of central atom $a_\ell$ with supp\,$a_\ell\subseteq B_\ell$, we have
\begin{equation*}
\begin{split}
I_1'&\le C\sum_{k\in\mathbb Z}\big[w_1(B_k)\big]^{{\alpha p}/n}
\Bigg(\sum_{\ell=k-1}^{\infty}\big|\lambda_\ell\big|\big\|a_\ell\big\|_{L^q_{w_2}}\Bigg)^p\\
&\le C\sum_{k\in\mathbb Z}\big[w_1(B_k)\big]^{{\alpha p}/n}
\Bigg(\sum_{\ell=k-1}^{\infty}\big|\lambda_\ell\big|\big[w_1(B_\ell)\big]^{-\alpha/n}\Bigg)^p\\
&\le C\sum_{k\in\mathbb Z}\big[w_1(B_k)\big]^{{\alpha p}/n}
\Bigg(\sum_{\ell=k-1}^{\infty}\big|\lambda_\ell\big|^p\big[w_1(B_\ell)\big]^{-{\alpha p}/{2n}}\Bigg)
\Bigg(\sum_{\ell=k-1}^\infty\big[w_1(B_\ell)\big]^{-{\alpha p'}/{2n}}\Bigg)^{p/{p'}}.
\end{split}
\end{equation*}
When $\ell\ge k-1$ with $k\in\mathbb Z$, then $B_{k-1}\subseteq B_\ell$. Since $w_1\in A_1$, as before, there exists a number $r>1$ such that $w_1\in RH_r$. Setting $\delta={(r-1)}/r>0$. Thus, by Lemma 2.2, we can see that
\begin{equation*}
\begin{split}
\sum_{\ell=k-1}^\infty\big[w_1(B_\ell)\big]^{-{\alpha p'}/{2n}}&=\big[w_1(B_{k-1})\big]^{-{\alpha p'}/{2n}}
\sum_{\ell=k-1}^\infty\left(\frac{w_1(B_{k-1})}{w_1(B_\ell)}\right)^{{\alpha p'}/{2n}}\\
&\le C\cdot \big[w_1(B_{k-1})\big]^{-{\alpha p'}/{2n}}
\sum_{\ell=k-1}^\infty \left(2^{(k-1)-\ell}\right)^{{\alpha\delta p'}/{2}}\\
&\le C\cdot \big[w_1(B_{k-1})\big]^{-{\alpha p'}/{2n}}\sum_{\ell=0}^\infty 2^{-{\ell\alpha\delta p'}/2}\\
&\le C\cdot \big[w_1(B_{k-1})\big]^{-{\alpha p'}/{2n}}.
\end{split}
\end{equation*}
Similarly,
\begin{equation*}
\sum_{k=-\infty}^{\ell+1}\big[w_1(B_{k-1})\big]^{{\alpha p}/{2n}}\big[w_1(B_\ell)\big]^{-{\alpha p}/{2n}}\le C,
\end{equation*}
where $C>0$ is an absolute constant which is independent of $\ell\in\mathbb Z$. Summarizing the estimates derived above, we thus obtain
\begin{equation*}
\begin{split}
I_1'&\le C\sum_{k\in\mathbb Z}\big[w_1(B_{k-1})\big]^{{\alpha p}/{2n}}
\Bigg(\sum_{\ell=k-1}^{\infty}\big|\lambda_\ell\big|^p\big[w_1(B_\ell)\big]^{-{\alpha p}/{2n}}\Bigg)\\
&= C\sum_{\ell\in\mathbb Z}\big|\lambda_\ell\big|^p
\Bigg(\sum_{k=-\infty}^{\ell+1}\big[w_1(B_{k-1})\big]^{{\alpha p}/{2n}}
\big[w_1(B_\ell)\big]^{-{\alpha p}/{2n}}\Bigg)\\
&\le C\sum_{\ell\in\mathbb Z}\big|\lambda_\ell\big|^p\\
&\le C\big\|f\big\|^p_{H\dot K^{\alpha,p}_q(w_1,w_2)}.
\end{split}
\end{equation*}
We now turn our attention to the estimate of $I_2'$. Observe that when $\ell\le k-2$, that is, $k\ge\ell+2$, then it follows immediately from the pointwise inequality (3.5) that
\begin{equation*}
\begin{split}
I_2'&\le C\sum_{k\in\mathbb Z}\big[w_1(B_k)\big]^{{\alpha p}/n}\Bigg(\sum_{\ell=-\infty}^{k-2}
\big|\lambda_\ell\big|\cdot
\frac{2^{\ell(n+\beta)}}{2^{k(n+\beta)}}\big[w_1(B_\ell)\big]^{-\alpha/n}\big[w_2(B_\ell)\big]^{-1/q}
\big[w_2(B_k)\big]^{1/q}\Bigg)^p\\
&=C\sum_{k\in\mathbb Z}\big[w_1(B_k)\big]^{{\alpha p}/n}\big[w_2(B_k)\big]^{p/q}\Bigg(\sum_{\ell=-\infty}^{k-2}\big|\lambda_\ell\big|\cdot
\frac{2^{\ell(n+\beta)}}{2^{k(n+\beta)}}\big[w_1(B_\ell)\big]^{-\alpha/n}\big[w_2(B_\ell)\big]^{-1/q}\Bigg)^p.
\end{split}
\end{equation*}
By using H\"older's inequality, we obtain that the above expression in the bracket is bounded by
\begin{equation*}
\begin{split}
&\Bigg(\sum_{\ell=-\infty}^{k-2}\big|\lambda_\ell\big|^p\cdot
\left(\frac{2^\ell}{2^k}\right)^{p(n+\beta)/2}\big[w_1(B_\ell)\big]^{-{\alpha p}/{2n}}\big[w_2(B_\ell)\big]^{-p/{2q}}\Bigg)\\
\times&\Bigg(\sum_{\ell=-\infty}^{k-2}\left(\frac{2^\ell}{2^k}\right)^{p'(n+\beta)/2}
\big[w_1(B_\ell)\big]^{-{\alpha p'}/{2n}}\big[w_2(B_\ell)\big]^{-p'/{2q}}\Bigg)^{p/{p'}}.
\end{split}
\end{equation*}
When $\ell\le k-2$ with $k\in \mathbb Z$, then we have $B_\ell\subseteq B_{k-2}\subseteq B_k$. Since $w_1,w_2\in A_1$, then it follows directly from Lemma 2.2 that
\begin{equation*}
\begin{split}
&\sum_{\ell=-\infty}^{k-2}\left(\frac{2^\ell}{2^k}\right)^{p'(n+\beta)/2}
\big[w_1(B_\ell)\big]^{-{\alpha p'}/{2n}}\big[w_2(B_\ell)\big]^{-{p'}/{2q}}\\
=&\,\big[w_1(B_{k})\big]^{-{\alpha p'}/{2n}}\big[w_2(B_{k})\big]^{-{p'}/{2q}}\\
&\times\sum_{\ell=-\infty}^{k-2}\left(\frac{2^\ell}{2^k}\right)^{p'(n+\beta)/2}
\left(\frac{w_1(B_{k})}{w_1(B_\ell)}\right)^{{\alpha p'}/{2n}}
\left(\frac{w_2(B_{k})}{w_2(B_\ell)}\right)^{p'/{2q}}\\
\le&\,C\cdot\big[w_1(B_{k})\big]^{-{\alpha p'}/{2n}}\big[w_2(B_{k})\big]^{-{p'}/{2q}}\\
&\times\sum_{\ell=-\infty}^{k-2}\left(\frac{2^\ell}{2^k}\right)^{p'(n+\beta)/2}
\left(\frac{2^{kn}}{2^{\ell n}}\right)^{{\alpha p'}/{2n}}\left(\frac{2^{kn}}{2^{\ell n}}\right)^{p'/{2q}}\\
\le&\,C\cdot\big[w_1(B_{k})\big]^{-{\alpha p'}/{2n}}\big[w_2(B_{k})\big]^{-{p'}/{2q}}
\cdot\sum_{\ell=2}^\infty\left(\frac{1}{2^\ell}\right)^{p'(n+\beta)/2-{\alpha p'}/{2}-{p'n}/{2q}}\\
\le&\,C\cdot\big[w_1(B_{k})\big]^{-{\alpha p'}/{2n}}\big[w_2(B_{k})\big]^{-{p'}/{2q}},
\end{split}
\end{equation*}
where the last inequality holds under our assumption that $\alpha<n(1-1/q)+\beta$. Similarly,
\begin{equation*}
\sum_{k=\ell+2}^\infty\left(\frac{2^\ell}{2^k}\right)^{p(n+\beta)/2}
\left(\frac{w_1(B_{k})}{w_1(B_\ell)}\right)^{{\alpha p}/{2n}}\left(\frac{w_2(B_{k})}{w_2(B_\ell)}\right)^{p/{2q}}\le C,
\end{equation*}
where $C>0$ is an absolute constant which is independent of $\ell\in\mathbb Z$. Summarizing the estimates derived above, we finally obtain
\begin{equation*}
\begin{split}
I_2'&\le C\sum_{k\in\mathbb Z}\big[w_1(B_{k})\big]^{{\alpha p}/{2n}}\big[w_2(B_{k})\big]^{p/{2q}}\\
&\times\Bigg(\sum_{\ell=-\infty}^{k-2}\big|\lambda_\ell\big|^p\cdot
\left(\frac{2^\ell}{2^k}\right)^{p(n+\beta)/2}\big[w_1(B_\ell)\big]^{-{\alpha p}/{2n}}\big[w_2(B_\ell)\big]^{-p/{2q}}\Bigg)\\
&\le C\sum_{\ell\in\mathbb Z}\big|\lambda_\ell\big|^p
\Bigg[\sum_{k=\ell+2}^\infty\left(\frac{2^\ell}{2^k}\right)^{p(n+\beta)/2}
\left(\frac{w_1(B_{k})}{w_1(B_\ell)}\right)^{{\alpha p}/{2n}}\left(\frac{w_2(B_{k})}{w_2(B_\ell)}\right)^{p/{2q}}\Bigg]\\
&\le C\sum_{\ell\in\mathbb Z}\big|\lambda_\ell\big|^p\\
&\le C\big\|f\big\|^p_{H\dot K^{\alpha,p}_q(w_1,w_2)}.
\end{split}
\end{equation*}
Therefore, summing up the above estimates for $I'_1$ and $I'_2$, we get the desired result. This completes the proof of Theorem 1.1.
\end{proof}

\begin{proof}[Proof of Theorem 1.2]
First we note that our assumptions $\alpha=n(1-1/q)+\beta$ and $0<\beta<1$ imply that $N=[\alpha+n(1/q-1)]=[\beta]=0$. According to Theorem 2.7, for every $f\in H\dot K^{\alpha,p}_q(w_1,w_2)$, we have the decomposition $f=\sum_{\ell\in\mathbb Z}\lambda_\ell a_\ell$, where $\sum_{\ell\in\mathbb Z}|\lambda_\ell|^p<\infty$ and each $a_\ell$ is a central $(\alpha,q,0;w_1,w_2)$-atom with supp\,$a_\ell\subseteq B_\ell=B(0,2^\ell)$. Then for any given $\sigma>0$, we write
\begin{equation*}
\begin{split}
&\sigma^p\cdot\sum_{k\in\mathbb Z}\big[w_1(B_k)\big]^{{\alpha p}/n}
w_2\Big(\Big\{x\in C_k:|\mathcal S_\beta(f)(x)|>\sigma\Big\}\Big)^{p/q}\\
\le&\sigma^p\cdot\sum_{k\in\mathbb Z}\big[w_1(B_k)\big]^{{\alpha p}/n}
w_2\bigg(\bigg\{x\in C_k:\sum_{\ell=k-1}^\infty\big|\lambda_\ell\big|\big|\mathcal S_\beta(a_\ell)(x)\big|>\sigma/2\bigg\}\bigg)^{p/q}\\
&+\sigma^p\cdot\sum_{k\in\mathbb Z}\big[w_1(B_k)\big]^{{\alpha p}/n}
w_2\bigg(\bigg\{x\in C_k:\sum_{\ell=-\infty}^{k-2}\big|\lambda_\ell\big|\big|\mathcal S_\beta(a_\ell)(x)\big|>\sigma/2\bigg\}\bigg)^{p/q}\\
=&J_1+J_2.
\end{split}
\end{equation*}
Since $w_2\in A_1$, then we have $w_2\in A_q$ for any $1<q<\infty$. Note that $0<p\le1$. Applying Chebyshev's inequality and Theorem A, we get
\begin{equation*}
\begin{split}
J_1&\le2^p\sum_{k\in\mathbb Z}\big[w_1(B_k)\big]^{{\alpha p}/n}
\Bigg(\sum_{\ell=k-1}^\infty\big|\lambda_\ell\big|\big\|\mathcal S_\beta(a_\ell)\chi_k\big\|_{L^q_{w_2}}\Bigg)^p\\
&\le 2^p\sum_{k\in\mathbb Z}\big[w_1(B_k)\big]^{{\alpha p}/n}
\Bigg(\sum_{\ell=k-1}^\infty\big|\lambda_\ell\big|^p\big\|\mathcal S_\beta(a_\ell)\big\|^p_{L^q_{w_2}}\Bigg)\\
&\le C\sum_{k\in\mathbb Z}\big[w_1(B_k)\big]^{{\alpha p}/n}
\Bigg(\sum_{\ell=k-1}^\infty\big|\lambda_\ell\big|^p\big\|a_\ell\big\|^p_{L^q_{w_2}}\Bigg).
\end{split}
\end{equation*}
Changing the order of summation yields
\begin{equation*}
J_1\le C\sum_{\ell\in\mathbb Z}\big|\lambda_\ell\big|^p
\Bigg(\sum_{k=-\infty}^{\ell+1}\big[w_1(B_k)\big]^{{\alpha p}/n}\big\|a_\ell\big\|^p_{L^q_{w_2}}\Bigg).
\end{equation*}
Following along the same lines as that of Theorem 1.1, we can also show that the series in the bracket is convergent. Furthermore, it is bounded by an absolute constant which is independent of $\ell\in\mathbb Z$. Hence
\begin{equation*}
J_1\le C\sum_{\ell\in\mathbb Z}\big|\lambda_\ell\big|^p\le C\big\|f\big\|^p_{H\dot K^{\alpha,p}_q(w_1,w_2)}.
\end{equation*}
On the other hand, observe that when $\ell\le k-2$, then for any $x\in C_k=B_k\backslash B_{k-1}$, by the pointwise inequality (3.5), we deduce that
\begin{equation*}
\begin{split}
\big|\mathcal S_\beta(a_{\ell})(x)\big|&\le C\cdot2^{\ell(n+\beta)}\big[w_1(B_\ell)\big]^{-\alpha/n}\big[w_2(B_\ell)\big]^{-1/q}
\cdot\frac{1}{|x|^{n+\beta}}\\
&\le C\cdot\frac{2^{\ell(n+\beta)}}{2^{k(n+\beta)}}\big[w_1(B_\ell)\big]^{-\alpha/n}\big[w_2(B_\ell)\big]^{-1/q}.
\end{split}
\end{equation*}
Since $B_j\subseteq B_{k-2}\subseteq B_k$ and $w_1, w_2\in A_1$, then it follows from our assumption $\alpha=n(1-1/q)+\beta$ and the inequality (3.6) that
\begin{align}
\big|\mathcal S_\beta(a_{\ell})(x)\big|&\le C\cdot\big[w_1(B_{k})\big]^{-\alpha/n}\big[w_2(B_{k})\big]^{-1/q}
\left(\frac{2^\ell}{2^{k}}\right)^{n+\beta}
\left(\frac{2^{kn}}{2^{\ell n}}\right)^{\alpha/n}\left(\frac{2^{kn}}{2^{\ell n}}\right)^{1/q}\notag\\
&\le C\cdot \big[w_1(B_{k})\big]^{-\alpha/n}\big[w_2(B_{k})\big]^{-1/q}.
\end{align}
Set $A_k=\big[w_1(B_{k})\big]^{-\alpha/n}\big[w_2(B_{k})\big]^{-1/q}$. We will consider the following two cases. If $\big\{x\in C_k:\sum_{\ell=-\infty}^{k-2}|\lambda_\ell||\mathcal S_\beta(a_\ell)(x)|>\sigma/2\big\}=\O$, then the inequality
\begin{equation*}
J_2\le C\big\|f\big\|^p_{H\dot K^{\alpha,p}_q(w_1,w_2)}
\end{equation*}
holds trivially. Now we assume that $\big\{x\in C_k:\sum_{\ell=-\infty}^{k-2}|\lambda_\ell||\mathcal S_\beta(a_\ell)(x)|>\sigma/2\big\}\neq\O$, then by the above inequality (3.7) and the fact that $0<p\le1$, we have
\begin{equation*}
\begin{split}
\sigma&< C\cdot A_k\Bigg(\sum_{\ell\in\mathbb Z}\big|\lambda_\ell\big|\Bigg)\\
&\le C\cdot A_k\Bigg(\sum_{\ell\in\mathbb Z}\big|\lambda_\ell\big|^p\Bigg)^{1/p}\\
&\le C\cdot A_k\big\|f\big\|_{H\dot K^{\alpha,p}_q(w_1,w_2)}.
\end{split}
\end{equation*}
It is easy to verify that $\lim_{k\to\infty}A_k=0$. Then for any fixed $\sigma>0$, we are able to find a maximal positive integer $K_\sigma$ such that
\begin{equation*}
\sigma<C\cdot A_{K_\sigma}\big\|f\big\|_{H\dot K^{\alpha,p}_q(w_1,w_2)}.
\end{equation*}
From the above discussions, we know that $B_{k}\subseteq B_{K_\sigma}$.Furthermore, by using Lemma 2.2 again, we obtain
\begin{equation*}
\frac{w_i(B_{k})}{w_i(B_{K_\sigma})}\le C\cdot\left(\frac{|B_{k}|}{|B_{K_\sigma}|}\right)^{\delta_i}\quad \mbox{for}\;\, i=1 \;\,\mbox{and}\;\, 2,
\end{equation*}
where $\delta_i>0$, $i=1,2$. Therefore
\begin{equation*}
\begin{split}
J_2&\le\sigma^p\cdot\sum_{k=-\infty}^{K_\sigma}\big[w_1(B_k)\big]^{{\alpha p}/n}\big[w_2(B_k)\big]^{p/q}\notag\\
&\le C\big\|f\big\|^p_{H\dot K^{\alpha,p}_q(w_1,w_2)}\sum_{k=-\infty}^{K_\sigma}
\left(\frac{w_1(B_k)}{w_1(B_{K_\sigma})}\right)^{{\alpha p}/n}
\left(\frac{w_2(B_k)}{w_2(B_{K_\sigma})}\right)^{p/q}\notag\\
&\le C\big\|f\big\|^p_{H\dot K^{\alpha,p}_q(w_1,w_2)}
\sum_{k=-\infty}^{K_\sigma}\left(\frac{1}{2^{(K_\sigma-k)n}}\right)^{\alpha\delta_1p/n+\delta_2p/q}\notag\\
&\le C\big\|f\big\|^p_{H\dot K^{\alpha,p}_q(w_1,w_2)}.
\end{split}
\end{equation*}
Combining the above estimates for $J_1$ and $J_2$, and then taking the supremum over all $\sigma>0$, we complete the proof of Theorem 1.2.
\end{proof}

\section{Proofs of Theorems 1.3 and 1.4}

In this section, we first establish the following three estimates which will be used in the proofs of our main theorems.

\newtheorem{prop}[theorem]{Proposition}

\begin{prop}
Let $w\in A_1$ and $0<\beta\le1$. Then for any $j\in\mathbb Z_+$, we have
\begin{equation*}
\big\|\mathcal S_{\beta,2^j}(a)\big\|_{L^2_w}\le C\cdot2^{jn/2}\big\|\mathcal S_\beta(a)\big\|_{L^2_w}.
\end{equation*}
\end{prop}

\begin{proof}
Since $w\in A_1$, then by Lemma 2.1, we know that for any $(y,t)\in{\mathbb R}^{n+1}_+$,
\begin{equation*}
w\big(B(y,2^jt)\big)=w\big(2^jB(y,t)\big)\le C\cdot2^{jn}w\big(B(y,t)\big),\quad j=1,2,\ldots.
\end{equation*}
Therefore, for any $j\in\mathbb Z_+$ and $0<\beta\le1$, we have
\begin{equation*}
\begin{split}
\big\|\mathcal S_{\beta,2^j}(a)\big\|_{L^2_w}^2&=\int_{\mathbb R^n}\bigg(\iint_{{\mathbb R}^{n+1}_+}\Big(A_\beta(a)(y,t)\Big)^2\chi_{|x-y|<2^j t}\frac{dydt}{t^{n+1}}\bigg)w(x)\,dx\\
&=\iint_{{\mathbb R}^{n+1}_+}\bigg(\int_{|x-y|<2^j t}w(x)\,dx\bigg)
\Big(A_\beta(a)(y,t)\Big)^2\frac{dydt}{t^{n+1}}\\
&\le C\cdot2^{jn}\iint_{{\mathbb R}^{n+1}_+}\bigg(\int_{|x-y|<t}w(x)\,dx\bigg)\Big(A_\beta(a)(y,t)\Big)^2\frac{dydt}{t^{n+1}}\\
&=C\cdot 2^{jn}\big\|\mathcal S_\beta(a)\big\|^2_{L^2_w}.
\end{split}
\end{equation*}
Taking square-roots on both sides of the above inequality, we are done.
\end{proof}

\begin{prop}
Let $w\in A_1$, $0<\beta\le1$ and $2<q<\infty$. Then for any $j\in\mathbb Z_+$, we have
\begin{equation*}
\big\|\mathcal S_{\beta,2^j}(a)\big\|_{L^q_w}\le C\cdot2^{jn/2}\big\|\mathcal S_\beta(a)\big\|_{L^q_w}.
\end{equation*}
\end{prop}

\begin{proof}
For any $j\in\mathbb Z_+$ and $0<\beta\le1$, it is easy to see that
\begin{equation}
\big\|\mathcal S_{\beta,2^j}(a)\big\|^2_{L^q_w}=\big\|\mathcal S_{\beta,2^j}(a)^2\big\|_{L^{q/2}_w}.
\end{equation}
Since $q/2>1$, by the duality argument, we then have
\begin{align}
&\big\|\mathcal S_{\beta,2^j}(a)^2\big\|_{L^{q/2}_w}\notag\\
=&\sup_{\|b\|_{L_w^{(q/2)'}}\le1}\left|\int_{\mathbb R^n}\mathcal S_{\beta,2^j}(a)(x)^2b(x)w(x)\,dx\right|\notag\\
=&\sup_{\|b\|_{L_w^{(q/2)'}}\le1}\left|\int_{\mathbb R^n}\bigg(\iint_{{\mathbb R}^{n+1}_+}\Big(A_\beta(a)(y,t)\Big)^2\chi_{|x-y|<2^j t}\frac{dydt}{t^{n+1}}\bigg)b(x)w(x)\,dx\right|\notag\\
=&\sup_{\|b\|_{L_w^{(q/2)'}}\le1}\left|\iint_{{\mathbb R}^{n+1}_+}\bigg(\int_{|x-y|<2^jt}b(x)w(x)\,dx\bigg)\Big(A_\beta(a)(y,t)\Big)^2 \frac{dydt}{t^{n+1}}\right|.
\end{align}
For $w\in A_1$, we denote the weighted maximal operator by $M_w$; that is
\begin{equation*}
M_w(f)(x)=\underset{x\in B}{\sup}\frac{1}{w(B)}\int_B|f(y)|w(y)\,dy,
\end{equation*}
where the supremum is taken over all balls $B$ which contain $x$. Hence, by using Lemma 2.1, we can get
\begin{align}
\int_{|x-y|<2^jt}b(x)w(x)\,dx&\le C\cdot2^{jn}w\big(B(y,t)\big)\cdot\frac{1}{w(B(y,2^jt))}\int_{B(y,2^jt)}b(x)w(x)\,dx\notag\\
&\le C\cdot2^{jn}w\big(B(y,t)\big)\underset{x\in B(y,2^jt)}{\inf}M_w(b)(x)\notag\\
&\le C\cdot2^{jn}\int_{|x-y|<t}M_w(b)(x)w(x)\,dx.
\end{align}
Substituting the above inequality (4.3) into (4.2) and then using H\"older's inequality together with the $L^{(q/2)'}_w$ boundedness of $M_w$, we thus obtain
\begin{equation*}
\begin{split}
\big\|\mathcal S_{\beta,2^j}(a)^2\big\|_{L^{q/2}_w}&\le C\cdot2^{jn}
\sup_{\|b\|_{L_w^{(q/2)'}}\le1}\left|\int_{\mathbb R^n}\mathcal S_\beta(a)(x)^2M_w(b)(x)w(x)\,dx\right|\\
&\le C\cdot2^{jn}\big\|\mathcal S_\beta(a)^2\big\|_{L^{q/2}_w}\sup_{\|b\|_{L_w^{(q/2)'}}\le1}\big\|M_w(b)\big\|_{L^{(q/2)'}_w}\\
&\le C\cdot2^{jn}\big\|\mathcal S_\beta(a)^2\big\|_{L^{q/2}_w}\\
&= C\cdot2^{jn}\big\|\mathcal S_\beta(a)\big\|^2_{L^q_w}.
\end{split}
\end{equation*}
This estimate together with (4.1) implies the desired result.
\end{proof}

\begin{prop}
Let $w\in A_1$, $0<\beta\le1$ and $1<q<2$. Then for any $j\in\mathbb Z_+$, we have
\begin{equation*}
\big\|\mathcal S_{\beta,2^j}(a)\big\|_{L^q_w}\le C\cdot2^{jn/q}\big\|\mathcal S_\beta(a)\big\|_{L^q_w}.
\end{equation*}
\end{prop}

\begin{proof}
We will adopt the same method given in \cite{torchinsky} to deal with the weighted case. For any $j\in\mathbb Z_+$ and $0<\beta\le1$, set $\Omega_\lambda=\big\{x\in\mathbb R^n:\mathcal S_\beta(a)(x)>\lambda\big\}$ and $\Omega_{\lambda,j}
=\big\{x\in\mathbb R^n:\mathcal S_{\beta,2^j}(a)(x)>\lambda\big\}$. We also set
\begin{equation*}
\Omega^*_\lambda=\Big\{x\in\mathbb R^n:M_w(\chi_{\Omega_\lambda})(x)>\frac{1}{2^{(jn+1)}\cdot[w]_{A_1}}\Big\}.
\end{equation*}
Observe that $w\big(\Omega_{\lambda,j}\big)\le w\big(\Omega^*_\lambda\big)+w\big(\Omega_{\lambda,j}\cap(\mathbb R^n\backslash\Omega^*_\lambda)\big)$. Thus, for any $j\in\mathbb Z_+$,
\begin{equation*}
\begin{split}
\big\|\mathcal S_{\beta,2^j}(a)\big\|^q_{L^q_w}&=\int_0^\infty q\lambda^{q-1}
\cdot w\big(\Omega_{\lambda,j}\big)\,d\lambda\\
&\le\int_0^\infty q\lambda^{q-1}\cdot w\big(\Omega^*_\lambda\big)\,d\lambda+\int_0^\infty q\lambda^{q-1}\cdot w\big(\Omega_{\lambda,j}\cap(\mathbb R^n\backslash\Omega^*_\lambda)\big)\,d\lambda\\
&=\mbox{\upshape I+II}.
\end{split}
\end{equation*}
The weighted weak type estimate of $M_w$ implies
\begin{equation}
\mbox{\upshape I}\le C\cdot2^{jn}\int_0^\infty q\lambda^{q-1}\cdot w\big(\Omega_\lambda\big)\,d\lambda\le C\cdot2^{jn}\big\|\mathcal S_\beta(a)\big\|^q_{L^q_w}.
\end{equation}
To estimate II, we now claim that the following inequality holds.
\begin{equation}
\int_{\mathbb R^n\backslash\Omega^*_\lambda}\mathcal S_{\beta,2^j}(a)(x)^2w(x)\,dx\le C\cdot2^{jn}\int_{\mathbb R^n\backslash\Omega_\lambda}\mathcal S_{\beta}(a)(x)^2w(x)\,dx.
\end{equation}
Assuming this claim for the moment, then it follows from Chebyshev's inequality and the inequality (4.5) that
\begin{equation*}
\begin{split}
w\big(\Omega_{\lambda,j}\cap(\mathbb R^n\backslash\Omega^*_\lambda)\big)&\le\lambda^{-2}\int_{\Omega_{\lambda,j}\cap(\mathbb R^n\backslash\Omega^*_\lambda)}\mathcal S_{\beta,2^j}(a)(x)^2w(x)\,dx\\
&\le\lambda^{-2}\int_{\mathbb R^n\backslash\Omega^*_\lambda}\mathcal S_{\beta,2^j}(a)(x)^2w(x)\,dx\\
&\le C\cdot2^{jn}\lambda^{-2}\int_{\mathbb R^n\backslash\Omega_\lambda}\mathcal S_{\beta}(a)(x)^2w(x)\,dx.
\end{split}
\end{equation*}
Hence
\begin{equation*}
\mbox{\upshape II}\le C\cdot2^{jn}\int_0^\infty q\lambda^{q-1}\bigg(\lambda^{-2}\int_{\mathbb R^n\backslash\Omega_\lambda}\mathcal S_{\beta}(a)(x)^2w(x)\,dx\bigg)d\lambda.
\end{equation*}
Changing the order of integration yields
\begin{align}
\mbox{\upshape II}&\le C\cdot2^{jn}\int_{\mathbb R^n}\mathcal S_\beta(a)(x)^2\bigg(\int_{|\mathcal S_\beta(a)(x)|}^\infty q\lambda^{q-3}\,d\lambda\bigg)w(x)\,dx\notag\\
&\le C\cdot2^{jn}\cdot\frac{q}{2-q}\big\|\mathcal S_\beta(a)\big\|^q_{L^q_w}.
\end{align}
Combining the above estimate (4.6) with (4.4) and taking $q$-th roots on both sides, we are done. So it remains to prove the inequality (4.5). Set $\Gamma_{2^j}(\mathbb R^n\backslash\Omega^*_\lambda)=\underset{x\in\mathbb R^n\backslash\Omega^*_\lambda}{\bigcup}\Gamma_{2^j}(x)$ and
$\Gamma(\mathbb R^n\backslash\Omega_\lambda)=\underset{x\in\mathbb R^n\backslash\Omega_\lambda}{\bigcup}\Gamma(x).$
For each given $(y,t)\in\Gamma_{2^j}(\mathbb R^n\backslash\Omega^*_\lambda)$, then by Lemma 2.1, we have
\begin{equation*}
w\big(B\big(y,2^jt\big)\cap\big(\mathbb R^n\backslash\Omega_\lambda^*\big)\big)\le C\cdot2^{jn}w\big(B(y,t)\big).
\end{equation*}
It is not difficult to check that $w\big(B(y,t)\cap\Omega_\lambda\big)\le\frac{w(B(y,t))}{2}$ and $\Gamma_{2^j}(\mathbb R^n\backslash\Omega^*_\lambda)\subseteq\Gamma(\mathbb R^n\backslash\Omega_\lambda)$. In fact, for any $(y,t)\in\Gamma_{2^j}(\mathbb R^n\backslash\Omega^*_\lambda)$, there exists a point $x\in \mathbb R^n\backslash\Omega^*_\lambda$ so that $(y,t)\in\Gamma_{2^j}(x)$. Then by Lemma 2.1, we can deduce
\begin{equation*}
\begin{split}
w\big(B(y,t)\cap\Omega_\lambda\big)&\le w\big(B(y,2^jt)\cap\Omega_\lambda\big)\\
&= \int_{B(y,2^jt)}\chi_{\Omega_\lambda}(z)w(z)\,dz\\
&\le [w]_{A_{1}}\cdot2^{jn}w\big(B(y,t)\big)\cdot
\frac{1}{w(B(y,2^jt))}\int_{B(y,2^jt)}\chi_{\Omega_\lambda}(z)w(z)\,dz.
\end{split}
\end{equation*}
Notice that $x\in B(y,2^jt)\cap(\mathbb R^n\backslash\Omega^*_\lambda)$. So we have
\begin{equation*}
\begin{split}
w\big(B(y,t)\cap\Omega_\lambda\big)\le [w]_{A_{1}}\cdot2^{jn}w\big(B(y,t)\big)\cdot
M_w(\chi_{\Omega_\lambda})(x)\le \frac{w(B(y,t))}{2}.
\end{split}
\end{equation*}
Hence
\begin{equation*}
\begin{split}
w\big(B(y,t)\big)&=w\big(B(y,t)\cap\Omega_\lambda\big)+w\big(B(y,t)\cap(\mathbb R^n\backslash\Omega_\lambda)\big)\\
&\le \frac{w(B(y,t))}{2}+w\big(B(y,t)\cap(\mathbb R^n\backslash\Omega_\lambda)\big),
\end{split}
\end{equation*}
which is equivalent to
\begin{equation*}
w\big(B(y,t)\big)\le 2\cdot w\big(B(y,t)\cap(\mathbb R^n\backslash\Omega_\lambda)\big).
\end{equation*}
The above inequality implies in particular that there is a point $z\in B(y,t)\cap(\mathbb R^n\backslash\Omega_\lambda)\neq\emptyset$. In this case, we have $(y,t)\in\Gamma(z)$ with $z\in \mathbb R^n\backslash\Omega_\lambda$, which implies $\Gamma_{2^j}(\mathbb R^n\backslash\Omega^*_\lambda)\subseteq\Gamma(\mathbb R^n\backslash\Omega_\lambda)$. Thus we obtain
\begin{equation*}
w\big(B\big(y,2^jt\big)\cap\big(\mathbb R^n\backslash\Omega_\lambda^*\big)\big)\le C\cdot2^{jn}w\big(B(y,t)\cap(\mathbb R^n\backslash\Omega_\lambda)\big).
\end{equation*}
Therefore
\begin{equation*}
\begin{split}
&\int_{\mathbb R^n\backslash\Omega^*_\lambda}\mathcal S_{\beta,2^j}(a)(x)^2w(x)\,dx\\
=&\int_{\mathbb R^n\backslash\Omega^*_\lambda}\bigg(\iint_{\Gamma_{2^j}(x)}\Big(A_\beta(a)(y,t)\Big)^2\frac{dydt}{t^{n+1}}
\bigg)w(x)\,dx\\
\le&\iint_{\Gamma_{2^j}(\mathbb R^n\backslash\Omega^*_\lambda)}\bigg(\int_{B(y,2^jt)\cap(\mathbb R^n\backslash\Omega_\lambda^*)}w(x)\,dx\bigg)\Big(A_\beta(a)(y,t)\Big)^2\frac{dydt}{t^{n+1}}\\
\le&\,C\cdot2^{jn}\iint_{\Gamma(\mathbb R^n\backslash\Omega_\lambda)}\bigg(\int_{B(y,t)\cap(\mathbb R^n\backslash\Omega_{\lambda})}w(x)\,dx\bigg)\Big(A_\beta(a)(y,t)\Big)^2\frac{dydt}{t^{n+1}}\\
\le&\,C\cdot2^{jn}\int_{\mathbb R^n\backslash\Omega_\lambda}\mathcal S_{\beta}(a)(x)^2w(x)\,dx,
\end{split}
\end{equation*}
which is just what we want. This finishes the proof of Proposition 4.3.
\end{proof}

We are now in a position to give the proof of Theorem 1.3.

\begin{proof}[Proof of Theorem 1.3]
In view of Theorem 2.7, as in the proof of Theorem 1.1 for the case of $0<p\le1$, we only need to show that for any central $(\alpha,q,0;w_1,w_2)$-atom $a$ with supp\,$a\subseteq B_{\ell}=B(0,2^\ell)$, $\ell\in\mathbb Z$, there exists a constant $C>0$ independent of $a$ such that $\big\|\mathcal G^*_{\lambda,\beta}(a)\big\|_{\dot K^{\alpha,p}_q(w_1,w_2)}\le C$. As before, we write
\begin{equation*}
\begin{split}
\big\|\mathcal G^*_{\lambda,\beta}(a)\big\|^p_{\dot K^{\alpha,p}_q(w_1,w_2)}=&\sum_{k=-\infty}^{\ell+1}
\big[w_1(B_k)\big]^{{\alpha p}/n}\big\|\mathcal G^*_{\lambda,\beta}(a)\chi_k\big\|^p_{L^q_{w_2}}\\
&+\sum_{k=\ell+2}^{\infty}\big[w_1(B_k)\big]^{{\alpha p}/n}\big\|\mathcal G^*_{\lambda,\beta}(a)\chi_k\big\|^p_{L^q_{w_2}}\\
=&\,K_1+K_2.
\end{split}
\end{equation*}
First, from the definition of $\mathcal G^*_{\lambda,\beta}$, we readily see that
\begin{align}
\big|\mathcal G^*_{\lambda,\beta}(a)(x)\big|^2=&\iint_{\mathbb R^{n+1}_+}\left(\frac{t}{t+|x-y|}\right)^{\lambda n}\Big(A_\beta(a)(y,t)\Big)^2\frac{dydt}{t^{n+1}}\notag\\
=&\int_0^\infty\int_{|x-y|<t}\left(\frac{t}{t+|x-y|}\right)^{\lambda n}\Big(A_\beta(a)(y,t)\Big)^2\frac{dydt}{t^{n+1}}\notag\\
&+\sum_{j=1}^\infty\int_0^\infty\int_{2^{j-1}t\le|x-y|<2^jt}\left(\frac{t}{t+|x-y|}\right)^{\lambda n}\Big(A_\beta(a)(y,t)\Big)^2\frac{dydt}{t^{n+1}}\notag\\
\le& C\bigg[\mathcal S_\beta(a)(x)^2+\sum_{j=1}^\infty 2^{-j\lambda n}\mathcal S_{\beta,2^j}(a)(x)^2\bigg].
\end{align}
Since $\lambda>2>\max\{1,2/q\}$ and $w_2\in A_1$. Thus, by applying Propositions 4.1--4.3, Theorem A and the inequality (4.7), we obtain
\begin{align}
\big\|\mathcal G^*_{\lambda,\beta}(a)\big\|_{L^q_{w_2}}
&\le C\Bigg(\big\|\mathcal S_\beta(a)\big\|_{L^q_{w_2}}+\sum_{j=1}^\infty 2^{-\frac{j\lambda n}{2}}\big\|\mathcal S_{\beta,2^j}(a)\big\|_{L^q_{w_2}}\Bigg)\notag\\
&\le C\Bigg(\big\|\mathcal S_\beta(a)\big\|_{L^q_{w_2}}+\sum_{j=1}^\infty 2^{-\frac{j\lambda n}{2}}\cdot\big[2^{\frac{jn}{2}}+2^{\frac{jn}{q}}\big]\big\|\mathcal S_{\beta}(a)\big\|_{L^q_{w_2}}\Bigg)\notag\\
&\le C\big\|a\big\|_{L^q_{w_2}}\Bigg(1+\sum_{j=1}^\infty2^{-\frac{j\lambda n}{2}}
\cdot\big[2^{\frac{jn}{2}}+2^{\frac{jn}{q}}\big]\Bigg)\notag\\
&\le C\big\|a\big\|_{L^q_{w_2}}.
\end{align}
Hence, for the term $K_1$, it follows directly from the above inequality (4.8) that
\begin{equation*}
\begin{split}
K_1&\le \sum_{k=-\infty}^{\ell+1}\big[w_1(B_k)\big]^{{\alpha p}/n}
\big\|\mathcal G^*_{\lambda,\beta}(a)\big\|^p_{L^q_{w_2}}\\
\end{split}
\end{equation*}
\begin{equation*}
\begin{split}
&\le C\sum_{k=-\infty}^{\ell+1}\big[w_1(B_k)\big]^{{\alpha p}/n}\big\|a\big\|^p_{L^q_{w_2}}.
\end{split}
\end{equation*}
Following along the same lines as in Theorem 1.1, we can also prove that $K_1\le C$.
On the other hand, in the proof of Theorem 1.1, for any fixed $\ell$ with $\ell\le k-2$ and $x\in C_k=B_k\backslash B_{k-1}$, we have already proved
\begin{equation}
\big|\mathcal S_\beta(a)(x)\big|\le
C\cdot\bigg(2^{\ell(n+\beta)}\big[w_1(B_\ell)\big]^{-\alpha/n}
\big[w_2(B_\ell)\big]^{-1/q}\bigg)\cdot|x|^{-n-\beta}.
\end{equation}
We are now going to estimate $\big|\mathcal S_{\beta,2^j}(a)(x)\big|$ for $j=1,2,\ldots$. Observe that if $x\in C_k=B_k\backslash B_{k-1}$, $k\ge\ell+2$ and $z\in B_\ell$, then we have $|z|\le\frac12|x|$. We also note that supp\,$\varphi\subseteq\{x\in\mathbb R^n:|x|\le1\}$, then for any given $z\in B_\ell$, $(y,t)\in\Gamma_{2^j}(x)$ and $x\in C_k$ with $k\ge\ell+2$, by a simple calculation, we can see that
\begin{equation*}
t+2^jt>|x-y|+|y-z|\ge|x-z|\ge|x|-|z|\ge\frac{|x|}{2}.
\end{equation*}
For every $j\in\mathbb Z_+$ and for all $x\in B_k\backslash B_{k-1}$ with $k\ge\ell+2$, it then follows from the preceding inequality (3.2) that
\begin{align}
\left|\mathcal S_{\beta,2^j}(a)(x)\right|&\le C\bigg(2^{\ell(n+\beta)}\big[w_1(B_\ell)\big]^{-\alpha/n}
\big[w_2(B_\ell)\big]^{-1/q}\bigg)
\left(\int_{\frac{|x|}{2^{j+2}}}^\infty\int_{|y-x|<2^jt}\frac{dydt}{t^{2n+2\beta+n+1}}\right)^{1/2}\notag\\
&\le C\bigg(2^{\ell(n+\beta)}\big[w_1(B_\ell)\big]^{-\alpha/n}
\big[w_2(B_\ell)\big]^{-1/q}\bigg)\left(
\int_{\frac{|x|}{2^{j+2}}}^\infty2^{jn}\cdot\frac{dt}{t^{2n+2\beta+1}}\right)^{1/2}\notag\\
&\le C\bigg(2^{\ell(n+\beta)}\big[w_1(B_\ell)\big]^{-\alpha/n}
\big[w_2(B_\ell)\big]^{-1/q}\bigg)\cdot\frac{2^{\frac{j(3n+2\beta)}{2}}}{|x|^{n+\beta}}.
\end{align}
Consequently
\begin{equation*}
\begin{split}
\big\|\mathcal S_{\beta,2^j}(a)\chi_k\big\|_{L^q_{w_2}}\le& C\cdot 2^{\ell(n+\beta)}\big[w_1(B_\ell)\big]^{-\alpha/n}
\big[w_2(B_\ell)\big]^{-1/q}\cdot2^{\frac{j(3n+2\beta)}{2}}\\
&\times\bigg(\int_{2^{k-1}<|x|\le2^k}\frac{w_2(x)}{|x|^{(n+\beta)q}}\,dx\bigg)^{1/q}\\
\le&C\cdot2^{\frac{j(3n+2\beta)}{2}}\big[w_1(B_\ell)\big]^{-\alpha/n}
\left(\frac{2^{\ell(n+\beta)}}{2^{k(n+\beta)}}\right)\left(\frac{w_2(B_k)}{w_2(B_\ell)}\right)^{1/q}.
\end{split}
\end{equation*}
Hence
\begin{align}
&\sum_{j=1}^\infty2^{-\frac{j\lambda n}{2}}\big\|\mathcal S_{\beta,2^j}(a)\chi_k\big\|_{L^q_{w_2}}\notag\\
\le&\,C\cdot\big[w_1(B_\ell)\big]^{-\alpha/n}
\left(\frac{2^{\ell(n+\beta)}}{2^{k(n+\beta)}}\right)
\left(\frac{w_2(B_k)}{w_2(B_\ell)}\right)^{1/q}\sum_{j=1}^\infty2^{-\frac{j(\lambda n-3n-2\beta)}{2}}\notag\\
\le&\,C\cdot\big[w_1(B_\ell)\big]^{-\alpha/n}
\left(\frac{2^{\ell(n+\beta)}}{2^{k(n+\beta)}}\right)
\left(\frac{w_2(B_k)}{w_2(B_\ell)}\right)^{1/q},
\end{align}
where the last inequality follows from the assumption that $\lambda>3+(2\beta)/n$. Substituting the above inequality (4.11) into the term $K_2$ and using (4.7), we thus obtain
\begin{equation*}
\begin{split}
K_2\le& \,C\sum_{k=\ell+2}^{\infty}\big[w_1(B_k)\big]^{{\alpha p}/n}
\bigg\{\big\|\mathcal S_\beta(a)\chi_k\big\|^p_{L^q_{w_2}}+\bigg(\sum_{j=1}^\infty 2^{-\frac{j\lambda n}{2}}\big\|\mathcal S_{\beta,2^j}(a)\chi_k\big\|_{L^q_{w_2}}\bigg)^p\bigg\}\\
\le&\,C\sum_{k=\ell+2}^{\infty}\left(\frac{2^{\ell p(n+\beta)}}{2^{kp(n+\beta)}}\right)
\left(\frac{w_1(B_k)}{w_1(B_\ell)}\right)^{{\alpha p}/n}\left(\frac{w_2(B_k)}{w_2(B_\ell)}\right)^{p/q}.
\end{split}
\end{equation*}
The rest of the proof is exactly the same as that of Theorem 1.1, we can get $K_2\le C$. Therefore, we conclude the proof of Theorem 1.3 for the case $0<p\le1$ by combining the above estimates for $K_1$ and $K_2$. Finally, by using the same arguments as in Theorem 1.1, we can also obtain the desired results for the case of $1<p<\infty$. We leave the details to the reader.
\end{proof}

\begin{proof}[Proof of Theorem 1.4]
According to Theorem 2.7 again, for every $f\in H\dot K^{\alpha,p}_q(w_1,w_2)$, we have the decomposition $f=\sum_{\ell\in\mathbb Z}\lambda_\ell a_\ell$, where $\sum_{\ell\in\mathbb Z}|\lambda_\ell|^p<\infty$ and each $a_\ell$ is a central $(\alpha,q,0;w_1,w_2)$-atom with supp\,$a_\ell\subseteq B_\ell=B(0,2^\ell)$. Then for any fixed $\sigma>0$, as in the proof of Theorem 1.2, we write
\begin{equation*}
\begin{split}
&\sigma^p\cdot\sum_{k\in\mathbb Z}\big[w_1(B_k)\big]^{{\alpha p}/n}
w_2\Big(\Big\{x\in C_k:\big|\mathcal G^*_{\lambda,\beta}(f)(x)\big|>\sigma\Big\}\Big)^{p/q}\\
\le&\sigma^p\cdot\sum_{k\in\mathbb Z}\big[w_1(B_k)\big]^{{\alpha p}/n}
w_2\bigg(\bigg\{x\in C_k:\sum_{\ell=k-1}^\infty\big|\lambda_\ell\big|\big|\mathcal G^*_{\lambda,\beta}(a_\ell)(x)\big|>\sigma/2\bigg\}\bigg)^{p/q}\\
&+\sigma^p\cdot\sum_{k\in\mathbb Z}\big[w_1(B_k)\big]^{{\alpha p}/n}
w_2\bigg(\bigg\{x\in C_k:\sum_{\ell=-\infty}^{k-2}\big|\lambda_\ell\big|\big|\mathcal G^*_{\lambda,\beta}(a_\ell)(x)\big|>\sigma/2\bigg\}\bigg)^{p/q}\\
=&K'_1+K'_2.
\end{split}
\end{equation*}
Note that $0<p\le1$ and $\lambda>2>\max\{1,2/q\}$. Applying Chebyshev's inequality and the inequality (4.8), we get
\begin{equation*}
\begin{split}
K'_1&\le 2^p\sum_{k\in\mathbb Z}\big[w_1(B_k)\big]^{{\alpha p}/n}
\Bigg(\sum_{\ell=k-1}^\infty\big|\lambda_\ell\big|\big\|\mathcal G^*_{\lambda,\beta}(a_\ell)\chi_k\big\|_{L^q_{w_2}}\Bigg)^p\\
&\le 2^p\sum_{k\in\mathbb Z}\big[w_1(B_k)\big]^{{\alpha p}/n}
\Bigg(\sum_{\ell=k-1}^\infty\big|\lambda_\ell\big|^p\big\|\mathcal G^*_{\lambda,\beta}(a_\ell)\big\|^p_{L^q_{w_2}}\Bigg)\\
&\le C\sum_{k\in\mathbb Z}\big[w_1(B_k)\big]^{{\alpha p}/n}
\Bigg(\sum_{\ell=k-1}^\infty\big|\lambda_\ell\big|^p\big\|a_\ell\big\|^p_{L^q_{w_2}}\Bigg).
\end{split}
\end{equation*}
Changing the order of summation gives us that
\begin{equation*}
K'_1\le C\sum_{\ell\in\mathbb Z}\big|\lambda_\ell\big|^p
\Bigg(\sum_{k=-\infty}^{\ell+1}\big[w_1(B_k)\big]^{{\alpha p}/n}\big\|a_\ell\big\|^p_{L^q_{w_2}}\Bigg).
\end{equation*}
Arguing as in the proof of Theorem 1.2, we can also show that
\begin{equation*}
K'_1\le C\big\|f\big\|^p_{H\dot K^{\alpha,p}_q(w_1,w_2)}.
\end{equation*}
We now turn to deal with $K'_2$. In this situation, it follows from the inequalities (4.7), (4.9) and (4.10) that
\begin{equation*}
\begin{split}
\big|\mathcal G^*_{\lambda,\beta}(a_\ell)(x)\big|&\le C\Bigg(\big|\mathcal S_\beta(a_\ell)(x)\big|+\sum_{j=1}^\infty 2^{-\frac{j\lambda n}{2}}\big|\mathcal S_{\beta,2^j}(a_\ell)(x)\big|\Bigg)\\
&\le C\bigg(2^{\ell(n+\beta)}\big[w_1(B_\ell)\big]^{-\alpha/n}
\big[w_2(B_\ell)\big]^{-1/q}\bigg)\big|x\big|^{-n-\beta}
\Bigg(1+\sum_{j=1}^\infty2^{-\frac{j(\lambda n-3n-2\beta)}{2}}\Bigg)\\
&\le C\bigg(2^{\ell(n+\beta)}\big[w_1(B_\ell)\big]^{-\alpha/n}
\big[w_2(B_\ell)\big]^{-1/q}\bigg)\big|x\big|^{-n-\beta},
\end{split}
\end{equation*}
where in the last inequality we have used the fact that $\lambda>3+(2\beta)/n$. Again, the rest of the proof is exactly the same as that of Theorem 1.2, we finally obtain
\begin{equation*}
K'_2\le C\big\|f\big\|^p_{H\dot K^{\alpha,p}_q(w_1,w_2)}.
\end{equation*}
Therefore, we conclude the proof of Theorem 1.4.
\end{proof}

\newtheorem*{rek}{Remark}
\begin{rek}
The corresponding results for non-homogeneous weighted Herz-type Hardy spaces can also be proved by atomic decomposition theory. The arguments are similar, so the details are omitted here.
\end{rek}

\end{document}